\documentclass{amsart}

\newtheorem{theorem}{Theorem}[section]

\newtheorem{proposition}[theorem]{Proposition}

\theoremstyle{definition}
\newtheorem{definition}[theorem]{Definition}
\newtheorem{example}[theorem]{Example}

\renewcommand{\P}{{\mathbb{P}}}
\newcommand{\Z}{{\mathbb{Z}}}
\newcommand{\N}{{\mathbb{N}}}
\newcommand{\Fc}{{\mathcal{F}}}
\newcommand{\Gc}{{\mathcal{G}}}
\newcommand{\Ic}{{\mathcal{I}}}
\newcommand{\Oc}{{\mathcal{O}}}
\newcommand{\Lc}{{\mathcal{L}}}

\newcommand{\reg}{{\mathrm{reg}}}

\newcommand{\unm}{\underline{m}}
\newcommand{\unn}{\underline{n}}
\newcommand{\und}{\underline{d}}
\newcommand{\unl}{\underline{l}}
\newcommand{\unp}{\underline{p}}

\newcommand\res[1]{{\lower1pt\hbox{$|$}}_{\raise.5pt\hbox{${\scriptstyle #1}$}}}

\numberwithin{equation}{section}

\begin{document}

\title{Regularity and Segre-Veronese Embeddings}

\author{David A.\ Cox} 
\address{Department of Mathematics and Computer Science, Amherst
College, Amherst, MA 01002-5000, USA} 
\email{dac@cs.amherst.edu}

\author{Evgeny Materov} 
\address{Department of Mathematics and Statistics, University of
Massachusetts, Amherst, MA 01003-9305, USA}
\email{materov@math.umass.edu}



\begin{abstract}
This paper studies the regularity of certain coherent sheaves that
arise naturally from Segre-Veronese embeddings of a product of
projective spaces.  We give an explicit formula for the regularity of
these sheaves and show that their regularity is subadditive.  We then
apply our results to study the Tate resolutions of these sheaves.
\end{abstract}

\date{\today}

\maketitle

\section{Regularity}

We will use the following generalized concept of regularity.

\begin{definition}
\label{Lregular}
Fix a projective variety $X$ and line bundles $L$, $B$ on $X$ such
that $B$ is generated by global sections.  Then a coherent sheaf $\Fc$
on $X$ is \emph{$L$-regular with respect to $B$} provided that
\[
  H^i(X,\Fc \otimes L \otimes B^{\otimes(-i)}) = 0 \quad \text{for
  all}\ i > 0. 
\]
Here $B^{\otimes(-i)}$ denotes the dual of $B^{\otimes i}$.  
\end{definition}

This notion of regularity is a special case the multigraded regularity
for sheaves introduced in \cite{hss}, which in turn is a modification
of the regularity defined in \cite{ms} for multigraded modules over a
polynomial ring.

We can relate Definition~\ref{Lregular} to the Castelnuovo-Mumford
regularity of a sheaf on projective space as follows.  Suppose that
$B$ is very ample and $\Lc = i_*\Fc$, where $i : X \to \P^N$ is the
projective embedding given by $B$.  Then, given an integer $p$, the
isomorphism
\[
H^i(\P^N,\Lc(p-i)) \simeq 
H^i(X,\Fc \otimes B^{\otimes p} \otimes B^{\otimes(-i)})
\]
shows that $\Fc$ is $B^{\otimes p}$-regular with respect to $B$ if and
only if $\Lc$ is $p$-regular as a coherent sheaf on $\P^N$.

Now fix $r$-tuples $\unl := (l_1,\dots,l_r)$ and $\und :=
(d_1,\dots,d_r)$ of positive integers.  These give the product variety
\[
X := \P^{l_1}\times\cdots\times\P^{l_r}
\]
of dimension $n := \sum_{i=1}^r l_i$ and the $\und$-uple
Segre-Veronese embedding 
\[
\nu_{\und} : X  \longrightarrow \P^N,\quad N =
  {\textstyle \prod_{k=1}^r \binom{l_k+d_k}{d_k}} - 1.
\]
The $r$-tuples $\unl$ and $\und$ will be fixed for the remainder of
the paper.

Given $\unm := (m_1,\ldots,m_r)$, consider the line bundle
\[
\Oc_X(\unm) = \Oc_X(m_1,\ldots,m_r)
\]
defined by
\[
  \Oc_X(\unm) := 
  p_1^*\Oc_{\P^{l_1}}(m_1)\otimes\cdots\otimes
  p_r^*\Oc_{\P^{l_r}}(m_r),
\]
where $p_i:\P^{l_1}\times\cdots\times\P^{l_r} \to \P^{l_i}$ is the
projection.  In this paper, we will study the regularity (in the sense
of Definition~\ref{Lregular}) of $\Oc_X(\unm)$ with respect to the
line bundles
\begin{align*}
L &:= \Oc_X(\unp) = \Oc_X(p_1,\dots,p_r), \ p_i \in \Z\\
B &:= \Oc_X(\und) = \Oc_X(d_1,\dots,d_r).
\end{align*}
For any nonempty subset $J \subseteq \{1,\ldots,r\}$ let $l_J$ denote
the sum $\sum_{j\in J}l_j$.  Here is our first main result.

\begin{theorem}
\label{theo_Lreg}
Let $L = \Oc_X(\unp)$ and $B = \Oc_X(\und)$ be as above.  Then
$\Oc_X(\unm)$ is $L$-regular with respect to $B$ if and only if
\begin{equation}
\label{pkineq}
  \max_{k\in J}\big\{p_k + m_k + l_k - l_J d_k\big\} \ge 0
\end{equation}
for all nonempty subsets $J \subseteq \{1,\ldots,r\}$.
\end{theorem}

\begin{proof}
Observe that
\begin{equation}
\label{rewriteLreg}
H^i(X,\Oc_X(\unm)\otimes L \otimes B^{\otimes(-i)}) = 
H^i(X,\Oc_X(\unm + \underline{p} -i\und)).
\end{equation}
and also that 
\[
  H^i(X,\Oc_X(\unm + \underline{p} -i\und)) = 
  \bigoplus_{i_1 + \cdots + i_r = i}
  \bigotimes_{k = 1}^r
  H^{i_k}(\P^{l_k},\Oc_{\P^{l_k}}(m_k + p_k-i d_k))
\]
by the K\"unneth formula.   Since
$H^{i_k}(\P^{l_k},\Oc_{\P^{l_k}}(j)) = 0$ when $i_k \ne 0, l_k$,
we may assume that $i = l_J$ for some $\emptyset \ne J
\subseteq \{1,\dots,r\}$ and that
\[
i_k = \begin{cases} l_k, & k \in J\\ 0, & k \notin J.\end{cases}
\]

First suppose that \eqref{pkineq} is satisfied for all $J \ne
\emptyset$.  Given such a $J$, pick $k \in J$ such that
\[
p_k + m_k + l_k - l_J d_k \ge 0.
\]
Then $m_k +p_k - l_J d_k \ge -l_k$, so that
$H^{l_k}(\P^{l_k},\Oc_{\P^{l_k}}(m_k +p_k -l_J d_k )) = 0$ by standard
vanishing theorems for line bundles on projective space.  By the above
analysis, it follows easily that \eqref{rewriteLreg} vanishes for $i >
0$.

Next suppose that 
\begin{equation}
\label{negative}
\max_{k \in J}\big\{p_k + m_k + l_k - l_J d_k\big\} < 0
\end{equation}
for some $J \ne \emptyset$.  Among all such subsets $J$, pick one of
maximum cardinality.  For this $J$, we
will prove that \eqref{rewriteLreg} is nonzero when $i = l_J$.  By the
K\"unneth formula, it suffices to show that
\begin{equation}
\label{toshow}
\begin{aligned}
H^{l_k}(\P^{l_k},\Oc_{\P^{l_k}}(m_k +p_k - l_J d_k )) &\ne 0, \quad k \in J\\
H^{\hskip1pt0\hskip1pt}(\P^{l_k},\Oc_{\P^{l_k}}(m_k +p_k - l_J d_k ))
&\ne 0, \quad k \notin J. 
\end{aligned}
\end{equation}
The first line of \eqref{toshow} is easy, since by Serre duality,
\[
H^{l_k}(\P^{l_k},\Oc_{\P^{l_k}}(m_k + p_k - l_J d_k )) \simeq 
H^{0}(\P^{l_k},\Oc_{\P^{l_k}}(- m_k - p_k +l_J d_k -l_k - 1))^*.
\]
This is nonzero provided $-m_k -p_k + l_J d_k -l_k - 1 \ge 0$, which for $k \in J$ follows immediately from 
\eqref{negative}.  For the
second line of \eqref{toshow}, suppose that $k \notin J$.  By the
maximality of $J$, we must have
\begin{equation}
\label{maxs}
\max_{s \in J\cup\{k\}}\big\{p_s + m_s + l_s - l_{J\cup\{k\}}
d_s\big\} \ge 0.   
\end{equation}
Note that $l_{J\cup\{k\}} = l_J + l_k$.  If the maximum in
\eqref{maxs} occurs at an element of $J$, say $s \in J$, then
\[
p_s + m_s + l_s - (l_J + l_k) d_s \ge 0,
\]
which is impossible since $s \in J$ and \eqref{negative} imply that
$p_s +m_s +l_s - I_J d_s < 0$.  Hence the maximum occurs at $k$, so
that
\[
p_k + m_k  + l_k - (l_J + l_k) d_k \ge 0.  
\]
This implies
\[
p_k + m_k  - l_J  d_k  \ge l_k d_k - l_k \ge 0.
\]
It follows
that the second line of \eqref{toshow} is nonzero, as desired.\end{proof}

If we fix the sheaf $\Oc_X(\unm)$ and let $L = \Oc_X(\unp)$ vary over
all $\unp \in \Z^r$, we get the \emph{regularity set}
\[
\mathrm{Reg}(\Oc_X(\unm)) = \{\unp \in \Z^r \mid \Oc_X(\unm)\
\text{is $\Oc_X(\unp)$-regular for $B = \Oc_X(\und)$}\}.
\]
This set is easy to describe.  Given a permutation $\sigma$ in the
symmetric group $S_r$, let $J(\sigma,k) \subseteq
\{1,\dots,r\}$ denote the subset
\[
J(\sigma,k) = \{\sigma(i) \mid i \ge \sigma^{-1}(k)\}, \quad 1 \le k \le r.
\]
Thus
\begin{align*}
J(\sigma,\sigma(1)) &= \{\sigma(1),\sigma(2),\dots,\sigma(r)\}\\
J(\sigma,\sigma(2)) &= \{\sigma(2),\dots,\sigma(r)\}\\
&\ \, \vdots\\
J(\sigma,\sigma(r)) &= \{\sigma(r)\}.
\end{align*}
Then define
\[
\unp{}_\sigma = (-m_1-l_1 + l_{J(\sigma,1)} d_1,\dots, -m_r-l_r +
l_{J(\sigma,r)} d_r) 
\]

\begin{proposition}
\label{regset}
The regularity set of $\Oc_X(\unm)$ is the union of $r!$ translates of
$\N^r$ given by
\[
\mathrm{Reg}(\Oc_X(\unm)) = \bigcup_{\sigma \in S_r} \big(\unp{}_\sigma
+ \N^r\big).
\]
\end{proposition}

\begin{proof}
First suppose that $\Oc_X(\unm)$ is $\Oc_X(\unp)$-regular for $B$.  We
build a permutation $\sigma \in S_r$ as follows.
Theorem~\ref{theo_Lreg} tells us that \eqref{pkineq} holds for all $J
\ne \emptyset$.  Hence, for $J = \{1,\dots,r\}$, there is $k_1$ such
that
\[
p_{k_1} + m_{k_1} + l_{k_1} - l_{J} d_{k_1} \ge 0.
\]
Setting $\sigma(1) = k_1$, this becomes 
\[
p_{k_1} \ge -m_{k_1}-l_{k_1} + l_{J(\sigma,k_1)} d_{k_1}.
\]
Applying \eqref{pkineq} to $J = \{1,\dots,r\}\setminus \{\sigma(1)\}$
gives $k_2 \ne k_1$ with
\[
p_{k_2} + m_{k_2} + l_{k_2} - l_{J} d_{k_2} \ge 0.
\]
Setting $\sigma(1) = k_2$, this becomes 
\[
p_{k_2} \ge -m_{k_2}-l_{k_2} + l_{J(\sigma,k_2)} d_{k_2}.
\]
Continuing in this way gives $\sigma \in S_r$ with $\unp \in \unp{}_\sigma
+ \N^r$.

Conversely, let $\unp \in \unp{}_\sigma + \N^r$ and take any nonempty
$J \subseteq \{1,\dots,r\}$.  Write $J = \{\sigma(i_1),\dots,
\sigma(i_t)\}$ for $i_1 < \cdots < i_t$.  Then set $k = \sigma(i_1)
\in J$ and observe that
\[
J \subseteq \{\sigma(i_1),\sigma(i_1+1),\dots, \sigma(r)\} =
J(\sigma,\sigma(i_1)) = J(\sigma,k).
\]
This implies that $l_J \le l_{J(\sigma,k)}$.  Since $\unp \in
\unp{}_\sigma + \N^r$, we obtain
\[
p_k \ge -m_{k}-l_{k} + l_{J(\sigma,k)} d_{k} \ge -m_{k}-l_{k} +
l_{J} d_{k},
\]
so that $p_k+m_k+l_k - l_J d_k \ge 0$.  This proves that \eqref{pkineq}
holds for $J$, and since $J \ne \emptyset$ was arbitrary, we have
$\unp \in \mathrm{Reg}(\Oc_X(\unm))$ by Theorem~\ref{theo_Lreg}.
\end{proof}

We next apply our results to study the Castelnuovo-Mumford regularity
of $\Oc_X(\unm)$ under the projective embedding given by the
Segre-Veronese map
\[
  \nu_{\und}: \P^{l_1}\times\cdots\times\P^{l_r} \longrightarrow \P^N
\]
coming from $\Oc_X(\und)$.   This gives the sheaf
$\Fc(\unm) := \nu_{\und*}\Oc_X(\unm)$ on $\P^N$.

\begin{theorem}
\label{theo_reg}
The sheaf $\Fc(\unm) = \nu_{\und*}\Oc_X(\unm)$ is $p$-regular on
$\P^N$ if and only if
\[
p \ge \max_{J \ne \emptyset}
  \left\{
  \min_{k\in J}\left\{l_J - \left\lfloor\frac{m_k + l_k}{d_k}
  \right\rfloor\right\}\right\}.
\]
Hence the regularity of $\Fc(\unm)$ is given by
\[
  \reg(\Fc(\unm)) = \max_{J \ne \emptyset}
  \left\{
  \min_{k\in J}\left\{l_J - \left\lfloor\frac{m_k + l_k}{d_k}
  \right\rfloor\right\}\right\}.
\]
\end{theorem}

\begin{proof}
Let $B = \Oc_X(\und)$.  As noted earlier, $\Fc(\unm)$ is $p$-regular
if and only if $\Oc_X(\unm)$ is $B^{\otimes p}$-regular with respect
to $B$.  Since $B^{\otimes p} = \Oc_X(p\und)$, Theorem~\ref{theo_Lreg}
implies that $\Fc(\unm)$ is $p$-regular
if and only if
\[
\max_{k \in J}\big\{p d_k + m_k + l_k - l_J d_k\big\} \ge 0
\]
for all subsets $J \ne \emptyset$ of $\{1,\dots,r\}$.  This is equivalent
to 
\[
p \ge \min_{k \in J}\left\{l_J  - \frac{m_k + l_k}{d_k}\right\} 
\]
for all $J \ne \emptyset$.  From here, the first assertion of the
theorem follows easily, and the second assertion is immediate.
\end{proof}

\begin{example}
When $r = 2$, the usual Segre embedding $\nu : \P^a \times \P^b \to
\P^{ab+a+b}$ gives $\Fc(k,l) = \nu_*\Oc_{\P^a \times \P^b}(k,l)$.
Since $\und = (1,1)$, Theorem~\ref{theo_reg} implies that
\begin{align*}
\reg(\Fc(k,l)) &= \max\{a-\lfloor k + a\rfloor, b-\lfloor l + b\rfloor,
\min\{ a+b - \lfloor k + a\rfloor, a+b - \lfloor l + b\rfloor\}\}\\
&= \max\{-k,-l,\min\{b-k,a-l\}\}\\ &=
\max\{-\min\{k,l\},\min\{b-k,a-l\}\}.
\end{align*}
This is the regularity formula from \cite[Lem.\ 3.1]{cm}.
\end{example}

\begin{example}
Let $\Ic$ be the ideal sheaf of $Y = \nu_{\und}(X) \subset \P^N$.
Following \cite{bs}, set $q_k = \lfloor \frac{l_k+1}{d_k}\rfloor$ and
$q_0 = \min_{1 \le k \le r}\{q_k\}$, and define
\[
\lambda = \begin{cases} n + 2 - q_0, & q _0 = q_k\ \text{for some $k$
    with}\ d_k | l_k+1\\ n + 1 - q_0, & \text{otherwise.}\end{cases}
\]
Recall that $n = \sum_{k=1}^r l_k = \dim(Y)$.  Then \cite[Lem.\
3.4]{bs} asserts that $\Ic$ is $\lambda$-regular.

We can derive this from Theorem~\ref{theo_reg} as follows.  We know
that $\Ic$ is $\lambda$-regular if and only if $\Oc_Y$ is $\lambda-1$
regular.  Note also that $\lambda$ can be defined more simply as
\[
\lambda = n+1-\min_{1\le k \le r} \left\{\left\lfloor
\frac{l_k}{d_k}\right\rfloor\right\}.
\]
If $J \ne \emptyset$, it follows easily that
\[
\lambda-1 \ge l_J - \min_{k\in J} \left\{\left\lfloor
\frac{l_k}{d_k}\right\rfloor\right\} \ge l_J - \max_{k\in J}
\left\{\left\lfloor \frac{l_k}{d_k}\right\rfloor\right\} = \min_{k\in
J} \left\{l_J - \left\lfloor \frac{l_k}{d_k}\right\rfloor\right\}.
\]
By Theorem~\ref{theo_reg}, we obtain $\lambda-1 \ge
\reg(\Fc(\underline{0})) = \reg(\Oc_Y)$, so that $\lambda \ge
\reg(\Ic)$, as claimed.

It is also easy to see that $\lambda - 1 > \reg(\Oc_Y)$ can occur.
For example, if $r = 2$ and $\und = (1,1)$, then one can show without
difficulty that
\[
\lambda-1 = \max\{l_1,l_2\} \ge \min\{l_1,l_2\} = \reg(\Oc_Y).
\]
\end{example}

\section{Subadditivity}

In this section we study how the regularity of $\Oc_X(\unm)$ and
$\Oc_X(\unm')$ compares to the regularity of the tensor product
\[
\Oc_X(\unm) \otimes \Oc_X(\unm') = \Oc_X(\unm + \unm').
\]
We first consider regularity as defined in
Definition~\ref{Lregular}. 

\begin{theorem}
\label{Lregadd}
Let $B = \Oc_X(\und)$.  If $\Oc_X(\unm)$ is $\Oc_X(\unp)$-regular for
$B$ and $\Oc_X(\unm')$ is $\Oc_X(\unp')$-regular for $B$, then
$\Oc_X(\unm + \unm')$ is $\Oc_X(\unp+\unp')$-regular for $B$.
\end{theorem}

\begin{proof}
We will use Theorem~\ref{theo_Lreg}.  Given a nonempty subset $J
\subseteq \{1,\dots,r\}$, it suffices to find $k \in J$ such that
\begin{equation}
\label{Ltoprove}
p_k + p'_k + m_k + m'_k + l_k - l_J d_k \ge 0.
\end{equation}
Using \eqref{pkineq} for this $J$ and the sheaves $\Oc_X(\unm)$ and
$\Oc_X(\unp)$, we know that
\[
\max_{s\in J}\big\{p_s + m_s + l_s - l_J d_s\big\} \ge 0. 
\]
Hence we can find $k \in J$ such that
\begin{equation}
\label{Lfirst}
p_k + m_k + l_k - l_J d_k \ge 0.
\end{equation}
Then using \eqref{pkineq} for $\{k\}$ and the sheaves $\Oc_X(\unm')$
and $\Oc_X(\unp')$, we also have
\[
p'_k + m'_k + l_k - l_k d_k \ge 0.
\]
Since $l_k \le l_k d_k$, this implies that
\begin{equation}
\label{Lsecond}
p'_k + m'_k \ge 0.
\end{equation}
Adding \eqref{Lfirst} and \eqref{Lsecond}, we obtain the desired
inequality \eqref{Ltoprove}.
\end{proof}

Theorem~\ref{Lregadd} gives the following result concerning the
Castelnuovo-Mumford regularity of the sheaves $\Fc(\unm) =
\nu_{\und*}\Oc_X(\unm)$ on projective space.

\begin{theorem}
\label{Fmreg}
Given any $\unm, \unm' \in \Z^r$, we have
\[
\reg(\Fc(\unm)) + \reg(\Fc(\unm')) \ge \reg(\Fc(\unm) \otimes \Fc(\unm')).
\]
\end{theorem}

\begin{proof}
It suffices to show that if $\Fc(\unm)$ is $p$-regular and
$\Fc(\unm')$ is $p'$-regular, then $\Fc(\unm+\unm')$ is
$(p+p')$-regular.  This follows immediately from Theorem~\ref{Lregadd}
and the already-noted equivalence
\[
\Fc(\unn)\ \text{is $q$-regular} \iff \Oc_X(\unn)\ \text{is
  $\Oc_X(q\und)$-regular for $B = \Oc_X(\und)$}.\qedhere 
\]
\end{proof}

In general, regularity is not subadditive, i.e., given coherent
sheaves $\Fc$ and $\Gc$ on $\P^N$, the inequality
\begin{equation}
\label{goodreg}
\reg(\Fc) + \reg(\Gc) \ge \reg(\Fc \otimes \Gc)
\end{equation}
may fail.  Here is an example due to Chardin.

\begin{example}
Let $R = k[x,y,z,t]$ and consider the ideals $I_n = \langle
z^n,t^n\rangle$ and $J_m = \langle x^{m-1}t-y^{m-1}z\rangle$.  When
$n,m \ge 3$, Example 1.13.6 of \cite{chardin1} implies that
\begin{equation}
\label{badreg}
\reg(I_n) + \reg(J_m) = m+2n-1 < \reg(I_n + J_m) = mn-1.
\end{equation}
As noted by Chardin \cite{chardin2}, this remains true when we work in
the larger ring $S = k[x,y,z,t,u,v]$.  The key point is that the
ideals $I_n$, $J_m$, and $I_n+J_m$ are saturated in $S$.

Now let $\Fc$ and $\Gc$ be the coherent sheaves associated to $S/I_n$
and $S/J_m$ respectively.  Then $\Fc \otimes \Gc$ is the sheaf
associated to $S/(I_n+J_m)$.  Since $I_n$ and $J_m$ are saturated in
$S$, \eqref{badreg} easily implies that
\[
\reg(\Fc) + \reg(\Gc) = m+2n-3 < \reg(\Fc \otimes \Gc) = mn-2
\]
when $n,m \ge 3$.  This shows that subadditivity fails in general.
\end{example}
 
However, there are certain situations where \eqref{goodreg} does hold,
such as when $\Fc$ or $\Gc$ is locally free (see \cite[Prop.\
1.8.9]{lazarsfeld}).  Theorem~\ref{Lregadd} shows that the sheaves
$\Fc(\unm)$ give another class of coherent sheaves for which
regularity is subadditive.

\section{Tate Resolutions}

By \cite{EFS}, a coherent sheaf $\Fc$ on the projective space $\P(W)$
has a Tate resolution
\[
\cdots \longrightarrow T^p(\Fc) \longrightarrow T^{p+1}(\Fc)
\longrightarrow \cdots
\]
of free graded $E$-modules, $E = \bigwedge W^*$, with terms 
\begin{equation}
\label{tate_form}
  T^p(\Fc) = {\textstyle \bigoplus_i}\, \widehat{E}(i-p)\otimes
  H^i(\Fc(p-i)).
\end{equation}
Here, $\widehat{E} = \mathrm{Hom}_K(E,K)$ is the dual over the base
field $K$.

Standard vanishing theorems imply
\[
T^p(\Fc) = \widehat{E}(-p)\otimes H^0(\Fc(p)) 
\quad \ \text{for}\ p \gg 0.
\]
Furthermore, if $\Fc = i_*\mathcal{E}$ for a locally free sheaf
$\mathcal{E}$ on an irreducible Cohen-Macaulay subvariety $Y
\hookrightarrow \P(W)$, then Serre duality and the same vanishing
theorems imply that
\[
T^p(\Fc) = \widehat{E}(n-p)\otimes H^n(\Fc(p-n)) \ \text{for}\ p \ll
0.
\]
In this situation, we define
\begin{align*}
p^+ &= \min\{p \mid T^p(\Fc) = \widehat{E}(-p)\otimes H^0(\Fc(p))\}\\
p^- &= \max\{p \mid T^p(\Fc) = \widehat{E}(n-p)\otimes H^n(\Fc(p-n))\}.
\end{align*}

In general, the differentials in the Tate resolution are hard to
describe.  The exceptions are the ``horizontal'' maps
\[
\widehat{E}(i-p)\otimes H^i(\Fc(p-i)) \longrightarrow 
\widehat{E}(i-p-1)\otimes H^i(\Fc(p+1-i)),
\]
which can be written down explicitly (see \cite {EFS}).  It follows
that the Tate resolution is known for $p \ge p^+$ and $p \le p^-$.  In
other words, the interesting part of the Tate resolution lies in the
range $p^- \le p \le p^+$.

For $X = \P^{l_1}\times \cdots \times \P^{l_r}$ and 
the Segre-Veronese embedding $\nu_{\und} : X \to \P^N = \P(W)$ coming from $\Oc_X(\und)$, we have the sheaf $\Fc(\unm)$, $\unm = (m_1,\dots,m_r)$.  When $r = 2$ and $\und = (1,1)$, the Tate resolution, including differentials, is described in \cite{cm}.  In the general case considered here, we restrict ourselves to giving explicit formulas for $p^+$ and $p^-$.  This determines the length of the interesting part of the Tate resolution.

To state our formulas, recall that $\unl = (l_1,\dots,l_r)$
and $\und = (d_1,\dots,d_r)$.  We also set $\underline{1} =
(1,\dots,1)$ and
\[
\underline{\widetilde{m}} = -\unm + n\und - \unl - \underline{1}
= (-m_1+nd_1-l_1-1,\dots,-m_r+nd_r-l_r-1).
\] 

\begin{theorem} For the sheaf $\Fc(\unm)$, we have
\begin{align*}
p^+ &= \phantom{-}\reg(\Fc(\unm)) = \phantom{-}\max_{J \ne \emptyset}
  \left\{
  \min_{k\in J}\left\{l_J - \left\lfloor\frac{{m}_k + l_k}{d_k}
  \right\rfloor\right\}\right\}\\
p^{-} &= -\reg(\Fc(\underline{\widetilde{m}})) = -\max_{J \ne \emptyset}
\left\{
\min_{k\in J}\left\{
\left\lceil
\frac{m_k + 1}{d_k}
\right\rceil - l_{J^c}
\right\}\right\},
\end{align*}
where $J^c$ denotes the complement of $J \subseteq \{1,\ldots,r\}$.
\end{theorem}

\begin{proof}
The equality $p^+ =\reg(\Fc(\unm))$ is immediate from
\eqref{tate_form} and the definition of regularity.  Then
Theorem~\ref{theo_reg} gives the desired formula for $p^+$.  Since
$\underline{\widetilde{m}} = -\unm + n\und - \unl -
\underline{1}$, Theorem~\ref{theo_reg} also implies that
\begin{align*}
\reg(\Fc(\underline{\widetilde{m}})) &= 
\max_{J \ne \emptyset}
  \left\{
  \min_{k\in J}\left\{l_J - \left\lfloor\frac{(-m_k+nd_k-l_k-1) + l_k}{d_k}
  \right\rfloor\right\}\right\}\\
&= \max_{J \ne \emptyset}
\left\{
\min_{k\in J}\left\{
\left\lceil
\frac{m_k + 1}{d_k}
\right\rceil - l_{J^c}
\right\}\right\},
\end{align*}
where the last equality follows easily using $l_J - n = -l_{J^c}$.

It remains to prove $p^- =
-\reg(\Fc(\underline{\widetilde{m}}))$.  Observe that $p \le
p^-$ if and only if 
\begin{equation}
\label{pminus}
  H^{n-i}(\Fc(\unm)(p - (n-i))) = 0 \quad \text{for all}\ i > 0.
\end{equation}
By Serre duality on $X = \P^{l_1}\times \cdots \times \P^{l_r}$,
\begin{align*}
H^{n-i}(\Fc(\unm)(p - (n-i))) &= H^{n-i}(X,\Oc_X(\unm + (p - (n-i)\und))\\
&= H^{i}(X, \Oc_X(-\unm -(p-n+i)\und - \unl -
\underline{1}))^{*}\\
&= H^{i}(\Fc(-\unm +(n-p-i)\und - \unl -
\underline{1}))^{*}\\ &=
H^{i}(\Fc(\underline{\widetilde{m}})(-p-i))^{*},
\end{align*}
so that \eqref{pminus} is equivalent to
\[
H^{i}(\Fc(\underline{\widetilde{m}})(-p -i)) = 0 \quad \text{for
  all}\ i > 0. 
\]
This vanishing is equivalent to $-p \ge
\reg(\Fc(\underline{\widetilde{m}}))$, so that
\[
p \le p^- \iff -p \ge
\reg(\Fc(\underline{\widetilde{m}})).
\]
Hence $p^- = -\reg(\Fc(\underline{\widetilde{m}}))$, as desired.
\end{proof}

\begin{example}
When $\und = (1,\dots,1)$ and $\unl = (l,\dots,l)$, we may assume that
$\unm = (m_1,\dots,m_r)$ with $m_1 \le \dots \le m_r$.  Then one can
show without difficulty that
\begin{align*}
p^+ &= \max_{1 \le i \le r}\{(i-1)l - m_i\}\\
p^- &= \min_{1 \le i \le r}\{(i-1)l - m_i\} -1.
\end{align*}
This makes it easy to compute $p^+ - p^-$, which is the length of the
interesting part of the Tate resolution.  For instance, if $m_1 =
\cdots = m_{r-1} = 0$ and $m_r = m$, where $m \ge (r-1)l$, then
\[
p^+ - p^- = m - l + 1,
\]
which can be arbitrarily large.  On the other hand, if $\unm$ is more
balanced, say $m_1 = \cdots = m_{r} = m$, then
\[
p^+ - p^- = (r-1) l + 1,
\]
independent of $m$.
\end{example}

\section*{Acknowledgements}

We are grateful to Jessica Sidman for suggesting that we recast the results of Sections 1 and~2 using the general
notion of regularity given in Definition~\ref{Lregular}.  Thanks also to Marc Chardin for his helpful comments about subadditivity
and regularity.  Finally, we are 
grateful to Nicolae Manolache for bringing the paper \cite{bs} to our
attention.

\end{document}